\documentclass[12pt, reqno]{amsart}
\usepackage{amssymb, mathrsfs}
\usepackage{latexsym}
\usepackage{amsmath}
\usepackage{amsthm}
\usepackage{amsfonts}
\usepackage{dsfont, upgreek}
\usepackage{mathtools}
\usepackage{epsfig}
\usepackage{amscd}
\usepackage{graphicx}
\usepackage{tikz-cd}
\usepackage{enumitem}
\setlist[enumerate]{itemsep=0.5ex}

\usepackage{color}
\usepackage{tikz}
\usetikzlibrary{patterns}

\usepackage[left=1.4in,right=1.4in,top=1.2in,bottom=1.2in]{geometry}

\usepackage{tikz}
\usetikzlibrary{decorations.markings}
\usetikzlibrary{arrows.meta}

\usepackage{extarrows}

\usepackage{adjustbox}
\usepackage{pgfplots}
\usepgfplotslibrary{colormaps}
\usepackage{slashed}

\usepackage[colorlinks=true,linkcolor=blue,citecolor=blue]{hyperref}

\usepackage[colorinlistoftodos,prependcaption,textsize=tiny]{todonotes}  

\usepackage[all,cmtip]{xy}
\theoremstyle{plain}
\newtheorem{theorem}{Theorem}[section]

\theoremstyle{definition}

\newtheorem*{claim*}{Claim}

\theoremstyle{remark}

\numberwithin{equation}{section}

\newcommand{\Sc}{\mathrm{Sc}}
\newcommand{\cl}{\mathbb{C}\ell}

\newcommand{\Bigwedge}{\mathord{\adjustbox{raise=.4ex, totalheight=.7\baselineskip}{$\bigwedge$}}}

\newcommand{\id}{\mathrm{id}}

\newcommand{\R}{\mathbb{R}}
\newcommand{\tr}{\mathrm{tr}}

\newcommand{\Z}{\mathbb{Z}}

\newcommand{\interior}[1]{%
	{\kern0pt#1}^{\mathrm{\,o}}%
}

\makeatletter
\let\save@mathaccent\mathaccent
\newcommand*\if@single[3]{%
	\setbox0\hbox{${\mathaccent"0362{#1}}^H$}%
	\setbox2\hbox{${\mathaccent"0362{\kern0pt#1}}^H$}%
	\ifdim\ht0=\ht2 #3\else #2\fi
}
\newcommand*\rel@kern[1]{\kern#1\dimexpr\macc@kerna}
\newcommand*\overbar[1]{\@ifnextchar^{{\wide@bar{#1}{0}}}{\wide@bar{#1}{1}}}
\newcommand*\wide@bar[2]{\if@single{#1}{\wide@bar@{#1}{#2}{1}}{\wide@bar@{#1}{#2}{2}}}
\newcommand*\wide@bar@[3]{%
	\begingroup
	\def\mathaccent##1##2{%
		\let\mathaccent\save@mathaccent
		\if#32 \let\macc@nucleus\first@char \fi
		\setbox\z@\hbox{$\macc@style{\macc@nucleus}_{}$}%
		\setbox\tw@\hbox{$\macc@style{\macc@nucleus}{}_{}$}%
		\dimen@\wd\tw@
		\advance\dimen@-\wd\z@
		\divide\dimen@ 3
		\@tempdima\wd\tw@
		\advance\@tempdima-\scriptspace
		\divide\@tempdima 10
		\advance\dimen@-\@tempdima
		\ifdim\dimen@>\z@ \dimen@0pt\fi
		\rel@kern{0.6}\kern-\dimen@
		\if#31
		\overline{\rel@kern{-0.6}\kern\dimen@\macc@nucleus\rel@kern{0.4}\kern\dimen@}%
		\advance\dimen@0.4\dimexpr\macc@kerna
		\let\final@kern#2%
		\ifdim\dimen@<\z@ \let\final@kern1\fi
		\if\final@kern1 \kern-\dimen@\fi
		\else
		\overline{\rel@kern{-0.6}\kern\dimen@#1}%
		\fi
	}%
	\macc@depth\@ne
	\let\math@bgroup\@empty \let\math@egroup\macc@set@skewchar
	\mathsurround\z@ \frozen@everymath{\mathgroup\macc@group\relax}%
	\macc@set@skewchar\relax
	\let\mathaccentV\macc@nested@a
	\if#31
	\macc@nested@a\relax111{#1}%
	\else
	\def\gobble@till@marker##1\endmarker{}%
	\futurelet\first@char\gobble@till@marker#1\endmarker
	\ifcat\noexpand\first@char A\else
	\def\first@char{}%
	\fi
	\macc@nested@a\relax111{\first@char}%
	\fi
	\endgroup
}
\makeatother

\begin{document}

\title{Rigidity of strictly convex domains in Euclidean spaces}
\author{Jinmin Wang}
\address[Jinmin Wang]{Department of Mathematics, Texas A\&M University}
\email{jinmin@tamu.edu}
\thanks{The first author is partially supported by NSF 1952693.}
\author{Zhizhang Xie}
\address[Zhizhang Xie]{ Department of Mathematics, Texas A\&M University }
\email{xie@math.tamu.edu}
\thanks{The second author is partially supported by NSF  1952693.}

\begin{abstract}
	In this paper, we prove a rigidity theorem for smooth strictly convex domains  in Euclidean spaces. 
\end{abstract}
\maketitle
	\section{Introduction}
	In  \cite{Wang:2022vf}, the authors proved a rigidity theorem for geodesic balls in Euclidean spaces \cite[Theorem 1.7]{Wang:2022vf}. The main purpose of this paper is to generalize such a rigidity theorem  to all strictly convex domains with smooth boundary  in $\mathbb R^n$. More precisely, we  have the following main theorem of the paper. 
\begin{theorem}\label{thm:balliso}
	Let $(M,g)$ be a strictly convex domain  with smooth boundary in $\R^n$ \textup{($n\geq 2$)}. Let $(N,\overbar g)$ be a spin Riemannian manifold with boundary and $f\colon N\to M$ be a spin map. If
		\begin{enumerate}[label=$(\arabic*)$]
		\item 
		$\Sc(\overbar{g})_x \geq \Sc(g)_{f(x)} = 0$ for all $x\in N$, 
		\item $H_{\overbar{g}}(\partial N)_y \geq  H_{g}(\partial M)_{f(y)}$ for all $y\in \partial N$, 
		\item $f$ is distance-non-increasing on $N$, 
		\item the degree of $f$ is nonzero,	
	\end{enumerate} 
then $f$ is an isometry. 
\end{theorem}

As a special case of   the above theorem, we see that,  given a strictly convex domain with smooth boundary in $\mathbb R^n$, one cannot  increase its metric,  its scalar curvature and  the mean curvature of its boundary simultaneously. 

The authors would like to thank Bo Zhu for helpful comments. 

\section{Proof of Theorem \ref{thm:balliso}}

The odd dimensional case reduces to the even dimensional case by  considering $f\times \id \colon N\times[0,1]\to M\times[0,1]$ and applying the index theory for manifolds with corners developed in  \cite{Wang:2021tq}. Hence without loss of generality,  we may assume that $N$ and $M$ are even-dimensional. 

 Let $S_N$ and $S_M$ be the  spinor bundles over $N$ and $M$. Consider the vector bundle $S_N\otimes f^* S_M$ over $N$, which carries a natural Hermitian metric and a unitary connection $\nabla$ compatible with the Clifford multiplication by elements of $\cl (TN)\otimes f^\ast \cl(TM)$. We will denote the Clifford multiplication of a vector $\overbar v\in TN$ by $\overbar c(\overbar v)$ and the Clifford multiplication of a vector $v\in f^*TM$ by $ c(v)$.

Consider the associated Dirac operator $D$ on $S_N\otimes f^* S_M$: 
$$D=\sum_{i=1}^n \overbar c(\overbar e_i)\nabla_{\overbar e_i}, $$
where $\{e_i\}$ is a local orthonormal basis of $TN$ and 
$$\nabla=\nabla^{S_N}\otimes 1+1\otimes f^*(\nabla^{S_M})$$
is the connection induced by the spinor connections on  $S_N$ and $f^\ast S_M$. 

Let $B$ be the local boundary condition on  $\partial N$ given by
\begin{equation}
	(\overbar \epsilon\otimes \epsilon)(\overbar c(\overbar e_n(x))\otimes c(e_n(x)))\varphi(x)=-\varphi(x),~\forall x\in\partial N,
\end{equation}
for all smooth sections $\varphi$ of $S_N\otimes f^* S_M$, where $\overbar \epsilon$ (resp. $\epsilon$) is the $\Z_2$-grading operator on $S_N$ (resp. $f^\ast S_M$), and $\overbar e_n(x)$ (resp. $e_n(x)$) is the unit inner normal vector at $x$ (resp. at $f(x)$).

\begin{proof}[{Proof of Theorem \ref{thm:balliso}}]

In even dimensions, it follows from  \cite[Theorem 1.1]{Lottboundary} that 
\begin{enumerate}[label=$(\arabic*)$]
	\item 
	$\Sc(\overbar{g})_x = \Sc(g)_{f(x)} = 0$ for all $x\in N$, 
	\item $H_{\overbar{g}}(\partial N)_y =  H_{g}(\partial M)_{f(y)}$ for all $y\in \partial N$. 
\end{enumerate} 
In fact, the proof of  \cite[Theorem 1.1]{Lottboundary} shows that 
there exists a non-zero parallel section $\varphi$ of $S_N\otimes f^\ast S_M$  satisfying the boundary condition $B$, i.e.,  $\nabla \varphi = 0$. As mentioned above, the odd dimensional case can be reduced to the even dimensional case by applying   \cite[Theorem 1.8] {Wang:2021tq}.  Moreover, since $\partial M$ is strictly convex, it follows from \cite[Proposition 3.1]{Wang:2022vf} that  $f\colon \partial N\to \partial M$ is a local  isometry. 

By \cite[Theorem 2.3]{Wang:2022vf}, $(N,\overbar g)$ is also flat. Since we shall need part of the argument of \cite[Theorem 2.3]{Wang:2022vf} for the proof of the current theorem, let us  repeat some key steps from the proof of \cite[Theorem 2.3]{Wang:2022vf} for the convenience of the reader.  Since $\nabla$ is a Hermitian connection and $\varphi$ is a nonzero parallel section, we may assume that $|\varphi|=1$ everywhere on $N$.
Since $(M, g)$ is  a strictly convex domain in $\R^n$, let $\{v_1,v_2,\ldots,v_n \}$ be the standard basis of $\mathbb R^n$, which we shall also view as parallel sections  of $TM$. As $M$ is strictly convex, there exist $y_1,y_2,\cdots,y_n$ in $\partial M$  such that the inner normal vector of $\partial M$ at $y_i$ is equal to $v_i$. Let $x_i\in \partial N$ be a point in the preimage $f^{-1}(y_i)$. 

Let $\Lambda$ be the collection of all subsets of $\{1,2,\ldots,n\}$. For $\lambda\in \Lambda$, we define
$$w_\lambda=\wedge_{i\in \lambda}v_i\in \Bigwedge^\ast TM$$
Note that $\{w_\lambda \}_{\lambda\in\Lambda}$ are parallel sections of $\Bigwedge^\ast TM$ that form an orthonormal basis of $\Bigwedge^\ast T_xM$ at every point $x\in M$.

With the section $\varphi$ of $S_N\otimes f^\ast S_M$  above, we define
$$\varphi_\lambda=(1\otimes c(w_\lambda))\varphi.$$
Since $w_\lambda$ is parallel, we see that $\varphi_\lambda$ is also parallel with respect to the connection $\nabla$. Note that $\nabla$ is a Hermitian connection that preserves the inner product on $S_N\otimes f^\ast S_M$. Therefore, for any pair of elements $\lambda, \mu\in \Lambda$, the function  $\langle\varphi_\lambda(x),\varphi_\mu(x)\rangle$ (as $x$ varies over $N$) is a constant function. 

\begin{claim*}
	The parallel sections $\{\varphi_\lambda\}_{\lambda\in\Lambda}$ are mutually orthogonal.
\end{claim*}

Note that $\dim(S_N\otimes f^\ast S_M)=2^n=|\Lambda|$, where $|\Lambda|$ is the cardinality of the set $\Lambda$.  Thus if the claim holds, then the curvature form of $S_N\otimes f^\ast S_M$ vanishes. Since $M$ is flat, the curvature form of $S_N\otimes f^\ast S_M$ is equal to the curvature form $R^{S_N}$ of $S_N$. By \cite[Theorem 2.7]{spinorialapproach}, we have
$$R^{S_N}_{X,Y}\sigma=\frac 1 2 R^{\overbar g}_{X,Y}\cdot \sigma,~ \textup{ for all } \sigma\in\Gamma(S_N) \textup{ and } X,Y\in\Gamma(TN),$$
where $R^{\overbar g}$ is the curvature form of the Levi--Civita connection on $TN$ with respect to $\overbar g$. It follows that $R^{\overbar g}=0$, that is, $\overbar g$ is flat.

Now we prove the claim. For each $\lambda \in \Lambda$ and $x\in M$, we denote by $V_\lambda$ the subspace in $T_xM \cong \mathbb R^n$ spanned by $\{v_i \}_{i\in\lambda}$.

Let $\lambda$ and $\mu$ be two distinct members of $\Lambda$.  Without loss of generality, we assume there exists $1\leq k \leq n$ such that $k\in \mu$ and $k\notin \lambda$. Equivalently, we have $v_k\in V_\lambda^\perp\cap V_\mu$. Let $\overbar v_k$ be the unit inner normal vector at $x_k$ of $\partial N$.  Recall that the section $\varphi$ satisfies the following boundary condition at $x_k$:
$$(\overbar \epsilon\otimes\epsilon)(\overbar c(\overbar v_k)\otimes c(v_k))\varphi(x_k)=-\varphi(x_k).$$
Note that for a vector $v\in\R^n$, we have
$$c(w_\lambda)c(v)=\begin{cases}
	(-1)^{|\lambda|}c(v) c(w_\lambda),& \textup{ if } v\in V_\lambda^\perp, \vspace{.3cm}  \\
	(-1)^{|\lambda|-1}c(v) c(w_\lambda),& \textup{ if } v \in V_\lambda,
\end{cases}$$
where $|\lambda|$ is the cardinality of the set $\lambda$. Therefore $\varphi_\lambda$ and $\varphi_\mu$ satisfy the following equations at $x_k$: 
$$(\overbar \epsilon\otimes\epsilon)(\overbar c(\overbar v_k)\otimes c(v_k))\varphi_\lambda(x_k)=-\varphi_\lambda(x_k), $$
$$(\overbar \epsilon\otimes\epsilon)(\overbar c(\overbar v_k)\otimes c(v_k))\varphi_\mu(x_k)=\varphi_\mu(x_k).$$
It follows that  $\langle\varphi_\lambda,\varphi_\mu\rangle$ vanishes at $x_k$, hence everywhere on $N$. This proves the claim, hence shows that $(N, \overbar g)$ is flat. 

We will  prove that $f$ preserves the second fundamental forms of $\partial N$ and $\partial M$. Before that, we shall need the following key observation. For any two points $x,y\in\partial N$, we denote by $\overbar  v_x$ and $\overbar  v_y$ the inner normal vectors at $x$ and $y$ in $\partial N$, and $ v_{x}$ and $ v_y$ the unit inner normal vectors at $f(x)$ and $f(y)$ in $\partial M$, respectively. Pick a path $\gamma$ connecting $x$ and $y$, then the parallel transport of $ v_y$ along $\gamma$ defines a vector field. In particular, $ v_y$ at $x$ is independent of the choice of the path, since $N$ is flat. We claim that
$$\langle \overbar  v_x,\overbar v_y\rangle=\langle  v_x, v_y\rangle.$$

Recall that the section $\varphi$ satisfies the boundary condition at $y$:
$$(\overbar \epsilon\overbar c(\overbar  v_y)\otimes 1)\varphi(y)=-(1\otimes \epsilon c( v_y))\varphi(y).$$
Since $\varphi$ is parallel and $ v_y$ is parallel along $\gamma$, we have that
$(\overbar \epsilon\overbar c(\overbar  v_y)\otimes 1)\varphi$ and $(1\otimes \epsilon c( v_y))\varphi$ are both parallel along $\gamma$. By the uniqueness of the parallel transport, we see that at $x$
$$(\overbar \epsilon\overbar c(\overbar  v_y)\otimes 1)\varphi(x)=-(1\otimes \epsilon c( v_y))\varphi(x).$$
The boundary condition at $x$ yields that
$$(\overbar \epsilon\overbar c(\overbar  v_x)\otimes 1)\varphi(x)=-(1\otimes \epsilon c( v_x))\varphi(x).$$
Therefore for any two real numbers $a_1$ and $a_2$, we have
$$(\overbar \epsilon\overbar c(a_1\overbar v_x+a_2\overbar  v_y)\otimes 1)\varphi(x)=-(1\otimes \epsilon c(a_1 v_x+a_2 v_y))\varphi(x).$$
Since $\varphi(x) \neq 0$, by taking the norm of both sides of the above equation, we see that $$|a_1\overbar v_x+a_2\overbar  v_y|=|a_1\overbar v_x+a_2\overbar  v_y|$$
for all $a_1, a_2\in \mathbb R$. It follows that $\langle \overbar  v_x,\overbar v_y\rangle=\langle  v_x, v_y\rangle$.

 	We now show that $f$ preserves the second fundamental forms of $\partial N$ and $\partial M$. For any $x\in\partial N$, we choose $n$ points $\{x_1,\ldots,x_n\}$ near $x$ in $\partial N$ so that the set $\{v_1, \ldots, v_n\}$ is linearly independent in $\R^n$, where $v_i$ is the unit inner normal vector of $\partial M$  at $f(x_i)$. Such a set of points  always exits since $\partial M$ is strictly convex. Let $\overbar v_i$ be the inner normal vector of $\partial N$ at $x_i$. The parallel transport of each $\overbar v_i$ gives a parallel vector field near $x$, which we still denote by $\overbar v_i$. The same argument above shows that 
 	$$\langle \overbar  v_i,\overbar v_j\rangle=\langle  v_i, v_j\rangle,~\forall i,j.$$
 	In particular, the set of vectors $\{\overbar v_1, \ldots, \overbar v_n\}$ is also linearly independent.
 	
 	Let $\overbar v$ (resp. $ v$)  be the unit inner normal vector field of $\partial N$ (resp. $\partial M$) near $x$ (resp. $f(x)$). The same argument  above again shows that  $\langle \overbar v,\overbar v_i\rangle=\langle v, v_i\rangle$ for all  $1\leq i\leq n$. Therefore, $\overbar v$ and $ v$ can be  written as linear combinations of $\overbar v_i$'s and $ v_i$'s with the same coefficients. In other words, there are smooth functions $k_1,\ldots,k_n$ defined on a neighborhood of $x$ in $\partial M$ such that
 	$$ v=\sum_{i=1}^n k_i v_i \text{ and }\overbar v=\sum_{i=1}^n (f^*k_i)\overbar v_i.$$
 	Let $\overbar w$ be an arbitrary  vector field tangent to $\partial N$. Since $ v_i$ and $\overbar  v_i$ are parallel, we have
 	$$\nabla_{f_*\overbar w}^M v=\sum_{i=1}^n f_*\overbar w(k_i)\cdot  v_i \textup{ and }\nabla^N_{\overbar w}\overbar v=\sum_{i=1}^n \overbar w(f^*k_i)\cdot \overbar v_i.$$
 	Note that $\overbar w(f^*k_i)=f^*(f_*\overbar w(k_i))$ by the chain rule. Therefore for any vector fields $\overbar w, \overbar u$ tangent to  $\partial N$, we have
 	$$\langle\nabla^N_{\overbar w}\overbar v,\nabla^N_{\overbar u}\overbar  v\rangle=\langle \nabla^M_{f_*\overbar w} v,\nabla^M_{f_*\overbar u} v\rangle.$$
 	Let $\{\overbar w_1, \ldots, \overbar w_{n-1}\}$ be a local orthonormal basis of $T\partial N$ near $x$. As $f$ is a local  isometry from $\partial N$ to $\partial M$, $\{f_*\overbar w_1, \ldots, f_*\overbar w_{n-1}\}$ is also a local orthonormal basis of $T\partial M$ near $f(x)$. Let $A=(A_{jk})$ and $\overbar A=(\overbar A_{jk})$ be the second fundamental forms of $\partial M$ and $\partial N$, that is,  
 	$$\nabla^N_{\overbar w_j}\overbar v=-\sum_{k=1}^{n-1}\overbar A_{jk}\overbar w_k,\text{ and }\nabla^N_{f_*\overbar w_j} v=-\sum_{k=1}^{n-1}A_{jk}f_*\overbar w_k.$$
 	Since $M$ is strictly convex, we assume without loss of generality that $A$ is a diagonal matrix with positive diagonal entries. 
 	
 	By rewriting   $\langle\nabla^N_{\overbar w}\overbar v,\nabla^N_{\overbar u}\overbar  v\rangle=\langle \nabla^M_{f_*\overbar w} v,\nabla^M_{f_*\overbar u} v\rangle$ in terms of the above  matrix entries,  we obtain that 
 	\[ \overbar A^{\, 2}=A^2. \] Since $\overbar A$ is symmetric, we have that
 	$O=\overbar AA^{-1}$ is an orthogonal matrix. Note that
 	$$\tr(\overbar A)=\tr(OA)=\sum_{j=1}^{n-1} O_{jj}A_{jj}\leq \sum_{j=1}^{n-1} A_{jj}=\tr(A),$$
 	where the second equality is because  $A$ is diagonal, and the third inequality is because $|O_{jj}|\leq 1$, since $O$ is orthogonal. 
 	Recall that the mean curvature of $N$ and $M$ are equal, that is, $\tr(\overbar A)=\tr(A)$. This implies that  $O_{jj} =1$ for each $1\leq j\leq n-1$. Therefore, $O$ is the identity matrix. It follows that $f$ preserves the second fundamental forms.

 	Let $\widetilde N$ be the universal cover of $N$ with the lift metric.  	
 	Fix a point $\tilde x$ in $\widetilde N$ and an orthonormal frame at $\tilde x$. The parallel transport of this orthonormal frame at $\tilde x$ defines a set of global  orthonormal basis $\{\overbar e_1, \ldots, \overbar e_n\}$ of $T\widetilde N$, where each $\overbar e_i$ is parallel. For any $\tilde y\in \widetilde N$, choose a smooth path $\gamma$ connecting $\tilde x$ and $\tilde y$, and define 
 	\[ \tilde y_i \coloneqq \int_{\gamma} \langle \dot \gamma, \overbar e_i\rangle  \]
 	where $\dot \gamma$ is the tangent vector of $\gamma$. Since $\widetilde N$ is simply connected and flat, the above integral is independent of the choice  of $\gamma$ among all smooth curves connecting $\tilde x$ and $\tilde y$. These functions $\tilde y_i\colon \widetilde N\to\R$   together give rise to a map $p\colon \widetilde N\to \R^n$ that is locally isometric. 
 	
 	Our next step is to show that $p\colon \widetilde N \to p(\widetilde N)$ is a Riemannian covering map.  First we   show that $\partial N$ has only one connected component. Otherwise, fix a connected component $C$ of $\partial N$. The distance from $C$ to $\partial N-C$ is positive, and is attained by some geodesic $\gamma$ connecting $x\in C$ and $y\in C'$, where $C'$ is a connected component of $\partial N - C$. Since the length of $\gamma$ is the minimum among all curves connecting $C$ to $\partial N - C$, it follows that $\gamma$ is orthogonal to both $C$ and $C'$, and lies in the interior of $N$ except the two end points. Let $U$ be a small neighborhood of $\gamma$. 
 	Since $N$ is flat, $U$ embeds isometrically into $\R^n$. Such an embedding maps $\gamma$  to a  line segment. Now since both $C$ and $C'$ are strictly convex, any line segment from $C$ to $C'$ inside $U$ parallel to $\gamma$ shorten the distance. This contradicts  the minimality of the chosen geodesic $\gamma$ and proves the claim.

  The exact same argument above also shows that $\partial \widetilde N$ has only one connected component.	Therefore $p(\partial\widetilde N)$ is connected in $\R^n$, and has the same metric and second fundamental form as $\partial M$. Thus $p(\partial\widetilde N)$ is a subset in $\R^n$ that only differs from $\partial M$ by an affine isometry. Indeed, this follows from the uniqueness of solutions to  the partial differential equations describing  $p(\partial \widetilde N)$ and $\partial M$ in $\mathbb R^n$. Moreover,  
 	$p$ restricted to $\partial\widetilde N$ is   a covering map, which in fact is an isometry if $n = \dim N\geq 3$.
 	
 	We claim that if $p(x)\in p(\partial\widetilde N)$, then $x\in\partial\widetilde N$. In other words, the map $p$ will never map an interior point of $\widetilde N$ to $p(\partial \widetilde N)$.   Assume to the contrary that there exists $x$ in the interior of $\widetilde N$ such that $p(x)\in p(\partial\widetilde N)$. The distance from $x$ to $\partial \widetilde N$ is attained by a unique geodesic segment  $\gamma$ from $x$ to  a point $y\in \partial \widetilde N$. Note that $\gamma$ is orthogonal to $\partial\widetilde N$. As $p$ is a local isometry, $p(\gamma)$ is a non-trivial line segment in $\R^n$ from  $p(x)$ to $p(y)$, which is orthogonal to $p(\partial \widetilde N)$ at $p(y)$. Since $\partial\widetilde N$ is convex, the vector in $\mathbb R^n$ from $p(y)$ to $p(x)$  is pointing inward (with respect to $p(\partial \widetilde N)$). Therefor $p(\gamma)$ lies entirely in the  inside\footnote{We have already shown that $p(\partial \widetilde N)$ is strictly convex smooth compact hypersurface in $\mathbb R^n$. It follows that $p(\partial \widetilde N)$ separates $\mathbb R^n$ into two parts. That is, $\mathbb R^n - p(\partial \widetilde N)$ consists of two connected components, exactly one of which is compact. We call the compact connected component of $\mathbb R^n - p(\partial \widetilde N)$ the inside of $p(\partial \widetilde N)$, and  the noncompact connected component of $\mathbb R^n - p(\partial \widetilde N)$ the outside of $p(\partial \widetilde N)$. } of  $p(\partial\widetilde N)$. Let $\alpha
 	\colon [0, 1] \to p(\partial \widetilde N)$ be a smooth path in $p(\partial \widetilde N)$ with $\alpha(0)=p(y)$ and $\alpha(1)=p(x)$. Since $p$ is a covering map on $\partial\widetilde N$, $\alpha$ lifts uniquely to a path $\widetilde \alpha$ such that $\widetilde \alpha(0)=y$. As $p$ is a local isometry near $y$, there is a unique geodesic $\gamma_t$ connecting $y$ and  $\widetilde\alpha(t)$ for all sufficiently small $t\in [0, 1]$, which is mapped isometrically under the map $p$ to the line segment connecting $p(y)$ and $\alpha(t)$. Since $p$ is a local isometry everywhere, we can continue the construction of such geodesics $\gamma_t$ for all $t\in [0, 1]$ . In particular, $p(\gamma_1)$ coincides with $p(\gamma)$. By construction, $\gamma_1$ and $\gamma$ have the same length and point towards the same direction starting from $y$. It follows that   $x$  coincides with the other end point $\widetilde \alpha(1)$ of $\gamma_1$, which lies in $\partial \widetilde N$ by construction. This contradicts the assumption that $x$ lies in the interior of $\widetilde N$. This finishes the proof of the claim. Note that the same argument also proves that every point in the inside of $p(\partial\widetilde N)$ admits at least one preimage in $\widetilde N$.
 	
 	The interior $\widetilde N - \partial\widetilde N$  of $\widetilde N$ is connected and $p(\widetilde N - \partial\widetilde N)$ is disjoint from $p(\partial\widetilde N)$, so $p(\widetilde N - \partial\widetilde N)$ lies entirely in the inside of $p(\partial\widetilde N)$. To summarize, we see that $p(\widetilde N)$ is precisely the region enclosed by   the hypersurface $p(\partial\widetilde N)$ in $\mathbb R^n$. As $p(\partial\widetilde N)$ coincides with $\partial M$ up to an affine isometry, $p(\widetilde N)$ coincides with $M$ up to an affine isometry. Without loss of generality, we may assume that $p(\widetilde N)=M$. 
 	
 	Now we show that $p\colon \widetilde N \to p(\widetilde N) = M$ is a covering map. Indeed, if $z$ is a point in the interior of $M$, then its preimage $p^{-1}(z)$ consists of only interior points of $\widetilde N$. Let $\varepsilon$ be the distance from $z$ to $\partial M$. Then the $\varepsilon$-neighborhood of each point in $p^{-1}(z)$ is mapped  isometrically under the map $p$ to the $\varepsilon$-neighborhood of $z$. In particular,  the $(\varepsilon/2)$-neighborhoods of points in $p^{-1}(z)$ are disjoint in  $\widetilde N$. The same holds when $z$ lies in $\partial M$, as each point in its preimage lies in $\partial \widetilde N$. 
 	
 	As $M$ is simply connected,  $p$ has to be the trivial covering map, hence  an isometry. In particular, $p\colon \partial \widetilde N\to \partial M$ is a homeomorphism. Let $\pi\colon\widetilde N\to N$ be the corresponding covering map for the universal cover $\widetilde N$ of $N$. Note that $f\circ \pi=p$ on $\partial\widetilde N$. Therefore the restriction  $\pi$ on $\partial\widetilde N$ is injective. It follows that $\widetilde N=N$ and $\pi$ is the identity map.

 	Now the map $h\coloneqq f\circ p^{-1}\colon M \to \widetilde N = N \to M$ is distance non-increasing, and equal to the identity map when restricted to $\partial M$. To prove the theorem, it suffices to show that any such map has to be the identity map on $M$. Let $x_1$ and $x_2$ be two arbitrary points on $\partial M$. Since $M$ is strictly convex, there is a unique   line segment $\ell$ connecting $x_1$ and $x_2$ that  lies entirely in $M$. Then $h(\ell)$ is a curve in $M$ connecting  $x_1$ and $x_2$, hence its length is at least the length of $\ell$. Since $h$ is distance non-increasing, it follows that $h$ maps   $\ell$ to itself isometrically. Note that all such line segments cover the whole $M$. This completes the proof.

\end{proof}


\end{document}